\UseRawInputEncoding
\documentclass[12pt,draft,reqno]{amsart}

\usepackage[all]{xy}
\usepackage{amsthm,array,amssymb,amscd,amsfonts,latexsym}

\usepackage{enumerate}
\usepackage{mathrsfs}

\theoremstyle{plain}

\theoremstyle{definition}

\newtheorem{nothing*}[theorem]{}
\newtheorem{subnothing*}[sub]{}

\theoremstyle{remark}

\newcommand{\cc}{\raise .4pt \hbox{{$\scriptstyle{\bullet}$}}}

\begin{document}

\title[
Rational differential forms]{Rational differential forms\\ on
the variety of flexes of plane cubics}
\author[Vladimir L. Popov]{Vladimir L. Popov}
\address{Steklov Mathematical Institute,
Russian Academy of Sciences, Gub\-kina 8,
Moscow 119991, Russia}
\email{popovvl@mi-ras.ru}

\begin{abstract}
We prove that for every positive integer $d$, there are no nonzero
regular differential $d$-forms on every smooth irreducible projective algebraic
 variety birationally isomorphic to the variety of flexes of plane cubics.
\end{abstract}

\maketitle

Below we use the standard notation from  \cite{GH}, \cite{S}.

Consider a three-dimensional complex vector space $V$ and fix a basis
$x_1, x_2, x_3$ of the dual space $V^*$. In the space ${\rm S}^3(V^*)$ of degree three forms on the space  $V$, all monomials of the form $x_1^{i_1}x_2^{i_2}x_3^{i_3}$, where $i_1+i_2+i_3=3$, constitute (for their fixed ordering) a basis.\;Let $\{a_{i_1, i_2, i_3}\mid i_1+i_2+i_3=3\}$ be the dual to it basis of the dual space
$({\rm S}^3(V^*))^*$.

The sets of forms  $\{x_j\}$ and $\{a_{i_1, i_2, i_3}\}$ are the projective coordinate systems on the  projective spaces
${\mathbf P}^2:=PV$ and ${\mathbf P}^9:= P{\rm S}^3(V^*)$
associated with
 $V$ and ${\rm S}^3(V^*)$ respectively.
 Let
\begin{equation*}
\pi\colon {\rm S}^3(V^*)\setminus \{0\}\to {\mathbf P}^9
\end{equation*}
be the natural projection.
Consider the following forms on the product
 ${\mathbf P}^2\times {\mathbf P}^9$:
\begin{equation*}
F:=\sum_{i_1+i_2+ i_3=3} a_{i_1, i_2, i_3} x_1^{i_1}x_2^{i_2}x_3^{i_3},
\quad
H:=
{\rm det}\Big(\frac{\partial^2F}{\partial x_i\partial x_j}\Big).
\end{equation*}
The closed subset $X$ of ${\mathbf P}^2\times {\mathbf P}^9$, defined by the system of equations $F=H=0$ was explored in the papers \cite{Ha}, \cite{K0}, \cite{K}.\;It is irreducible,
nine-dimensional and has singularities
 (as a matter of fact, $X$ even is non-normal) \cite{K}.\;Let
$${\mathbf P}^2\xleftarrow{p_2} X\xrightarrow{p_9} {\mathbf P}^9$$
be the natural projections.\;For every irreducible form $f\in {\rm S}^3(V^*)$ such that the equation
$f=0$ defines on ${\mathbf P}^2$ a smooth cubic $C$, the fiber of the morphism
 $p_9$ over the point $\pi(f)\in {\mathbf P}^9$ is
 a set of nine points.\;The image of this set under the projection $p_2$
is exactly the set of all (nine) flexes  (inflection points) of the cubic $C$.\;Therefore, the set of points of $X$ in general position is identified with the set of all pairs
 $(C, c)$, where $C$ is a smooth cubic on ${\mathbf P}^2$, and $c$ is its flex.
 For this reason, $X$ is called (see.\;\cite{K}) {\it the variety of flexes
 of plane cubics}.

Let $Y$ be a smooth irreducible projective algebraic variety birational\-ly isomorphic to the algebraic variety $X$. The main result (Theorem 4)
of the paper \cite{K} is the claim that the irregularity of the variety $Y$ vanishes:
\begin{equation}\label{H1}
H^1(Y, \mathcal O_Y)=0.
\end{equation}

Below is proved a theorem, a special case of which  is the equality\;\eqref{H1}:
\vskip 2mm
\noindent{\sc Theorem.}
{\it Maintain the above notation.\;Then for every positive inte\-ger $d$,
the following properties hold:
\begin{enumerate}[\hskip 2.2mm\rm(i)]
\item $H^d(Y, \mathcal O_Y)=0$;
\item there are no nonzero regular differential $d$-forms on the variety\;$Y$\!.
\end{enumerate}
}

\noindent {\it Proof.} First, we note that properties (i) and (ii) are equivalent.
Indeed, according to the Hodge decomposition,
$H^q(Y, \Omega_Y^p)=\overline{H^p(Y, \Omega_Y^q)}$ for all integers $p, q\geqslant 0$
(see\;\cite{GH}).\;In view of $\Omega^0_Y={\mathcal O}_Y$, for
$q=d, p=0$, this gives
$H^d(Y, {\mathcal O}_Y)=\overline{H^0(Y,\Omega^d_Y)}$, whence it follows
that  property (i) is equivalent to the equality
$H^0(Y,\Omega^d_Y)=0$, which, in turn, is a reformulation of  property\;(ii).

We now prove that property (ii) holds.

The natural actions of the group ${\rm GL}(V)$ on $V$, $V^*$, and ${\rm S}^3(V^*)$  induce its actions on ${\mathbf P}^2$, ${\mathbf P}^9$ and ${\mathbf P}^2\times {\mathbf P}^9$, whose inefficiency kernels contain the group of scalar transformations  $Z:=\{\alpha\,{\rm Id}_V\mid \alpha\in {\mathbf C}^\times\}$. Therefore these latter actions induce the actions of the projective group  $G:={\rm PGL}(V)={\rm GL}(V)/Z$
on  ${\mathbf P}^2$, ${\mathbf P}^9$, and ${\mathbf P}^2\times {\mathbf P}^9$.

The variety $X$ is invariant under the specified action of the group  $G$ on ${\mathbf P}^2\times {\mathbf P}^9$.\;The classical Hesse's results
(\cite{H}; see also
\cite[pp.\;291--299]{BK}) yield the following statements:

(a) For every smooth cubic on ${\mathbf P}^2$, there is a transformation from $G$, which maps this cubic to the cubic on ${\mathbf P}^2$, defined by the equation \begin{equation}\label{He}
x_1^3+x_2^3+x_3^3+\lambda x_1x_2x_3=0,\quad \lambda\in \mathbf C.
\end{equation}

(b) For every $\lambda\in \mathbf C$, distinct from  $-3$, $-3\varepsilon$, $-3\varepsilon^{-2}$, where $\varepsilon=e^{2\pi i/3}$,  the cubic on  ${\mathbf P}^2$ defined by equation  \eqref{He}, is smooth and has exactly nine flexes, one of which is the point
$$q:=(0:-1: 1).$$

(c)  For every smooth cubic $C$ on ${\mathbf P}^2$ defined by equation
\eqref{He}, and every its flex $c$, there is an element   $g\in G$ such that $g(C)=C$ and $g(c)=q$.

We now consider the morphism $\varphi\colon {\mathbf A}^1\to {\mathbf P}^2\times {\mathbf P}^9$ defined by the formula
\begin{equation*}
\varphi(\lambda):=\big(q, \pi(x_1^3+x_2^3+x_3^3+\lambda x_1x_2x_3)\big) \quad \mbox{for every $\lambda\in{\mathbb A}^1={\mathbb C}$}.
\end{equation*}

It follows from (b) that $\varphi({\mathbf A}^1)$ lies in $X$, and from (a) and (c) that the morpism
\begin{equation}\label{G}
G\times {\mathbf A}^1\to X,\quad (g, \lambda)\mapsto g(\varphi(\lambda)),
\end{equation}
is dominant. Since the underlying variety of every connected affine algebraic group is rational (see \cite[Cor.\;2]{C}), it follows from the domi\-nance of morphism \eqref{G} that
the variety $X$, and hence $Y$, is unirational. Therefore there is a dominant rational map  $\psi\colon {\mathbf P}^9\dashrightarrow Y$.\;In view of the smoothness of ${\mathbf P}^9$ and $Y$ and the projectivity of $Y$, the induced homomorphism of the spaces of rational differential $d$-forms
$$\psi^*\colon \Omega^d(Y)\to \Omega^d({\mathbf P}^9)$$ defines an embedding
of the spaces of regular differential $d$-forms $$\Omega^d[Y]\hookrightarrow \Omega^d[{\mathbf P}^9]$$
(see\,\cite[Chap.\,III, \S6.1, Thm.\,2]{S}).\;Since
$\Omega^r[{\mathbf P}^n]=0$ for any positive $r$ and $n$ (see\;\cite[Chap.\;0, Sect.\;7]{GH}),
this implies that $\Omega^d[Y]=0$, i.e., that statement (ii) holds.\quad $\square$

\vskip 2mm

\noindent{\sc Remark.} Along the way of proof, unirationality of $X$ is proved.


\begin{thebibliography}{9}
\bibitem{BK} E.\;Brieskorn, H.\;Kn\"orrer, {\it Plane Algebraic Curves}, Birkh\"auser, Basel, 1986.
    \bibitem{C} C.\;Chevalley, {\it On algebraic group varieties}, J. Math. Soc. Japan {\bf 6} (1954), nos. 3--4, 303--324.
    \bibitem{GH} P.\;Griffiths, J.\;Harris, {\it Principles of Algebraic Geometry}, Wiley, New York, 1978.
 \bibitem{Ha}       J.\;Harris, {\it Galois groups of enumerative problems}, Duke Math. J.
 {\bf 46} (1979), no. 4, 685--724.
  \bibitem{H}      O.\;Hesse, {\it \"Uber die Elimination der Variabeln aus drei algebraischen Glei\-chun\-gen vom zweiten Grade mit zwei Variabeln}, J. Reine Angew. Math. {\bf 28} (1844), 68--96.
 \bibitem{K0}     V. S. Kulikov, {\it The Hesse curve of a Lefschetz pencil of plane curves},
 Russian Math. Surveys {\bf 72} (2017), no. 3, 574.
\bibitem{K} Vik.\;S.\;Kulikov, {\it On the variety of the inflection points of plane cubic curves}, {\tt arXiv:1810.01705v1} (3 Oct 2018).
  \bibitem{S}  I. R. Shafarevich,\;{\it Basic Algebraic Geometry}, Springer, Heidelberg, 2013, 2007.


\end{thebibliography}
\end{document}